\documentclass[10pt]{article}
\usepackage{amssymb}
\begin{document}
\newtheorem{theorem}{Theorem}[section]
\newtheorem{lemma}{Lemma}[section]
\newtheorem{remark}{Remark}[section]
\def\operatorname#1{{\rm#1\,}}
\def\text#1{{\hbox{#1}}}
\def\qedbox{\hbox{$\rlap{$\sqcap$}\sqcup$}}
\def\Cal#1{{\mathcal{#1}}}
%\input amstex
%\documentstyle{amsppt}
%\def\version{File S06v4a.tex last changed 05April2001 by PBG}
%\def\nmonth{\ifcase\month\ \or January\or
%   February\or March\or April\or May\or
%   June\or July\or August\or
%   September\or October\or November\else December\fi}
%\NoRunningHeads
%\NoBlackBoxes
%\font\small=cmr7\def\nmonth{\ifcase\month\ \or January\or
%February\or March\or April\or May\or June\or July\or August\or 
%September\or October\or November\else December\fi} 
%\newcount\qctqq\newcount\qctdat
%\def\date#1{{\global\qctdat=1}\footline=
%{\ifnum\number\qctdat =1{\global\qctdat=0
%\small\hfill\TeX\ \number
%\day\ \nmonth\ \number\year\ #1\ P. Gilkey}\fi\hfil}}
%\def\DATE{\nmonth\ \number\day\ \number\year}
%\newcount\qcts\newcount\qcta\qcta=
%  1\newcount\qcteq\newcount\qcthead
\def\setref#1{}
\def\sethead#1{}
\def\seteqn#1{}
\def\eqnpbg{\ifnum\qcteq=1 a\else\ifnum\qcteq=2 b\else
   \ifnum\qcteq=3 c\else\ifnum\qcteq=4 d\else
   \ifnum\qcteq=5 e\else\ifnum\qcteq=6 f\else\ifnum\qcteq=7 g\else
   \ifnum\qcteq=8 h\else\ifnum\qcteq=9 i\else
   \ifnum\qcteq=10 j\else\ifnum\qcteq=11 k\else\ifnum\qcteq=12 l\else
   \ifnum\qcteq=13 m\else\ifnum\qcteq=14 n\else\ifnum\qcteq=15 o\else
   \ifnum\qcteq=16 p\else\ifnum\qcteq=17 q\else
   \ifnum\qcteq=18 r\else\ifnum\qcteq=19 s\else\ifnum\qcteq=20 t\else
   \ifnum\qcteq=21u\else\ifnum\qcteq=22 v\else
   \ifnum\qcteq=23 w\else\ifnum\qcteq=24 x\else
   \ifnum\qcteq=25 y\else\ifnum\qcteq=26 z\else *
    \fi\fi\fi\fi\fi\fi\fi\fi\fi\fi\fi\fi\fi\fi\fi\fi\fi\fi\fi\fi\fi\fi\fi\fi\fi\fi}
%\newcount\qct\newcount\qcta\qct=0\qcta=1
%\def\pbgkey#1{\key{#1}\global\advance\qct by\number\qcta
%       \immediate\write3{\string\def\string#1\string{\number\qct\string}}}
\def\pbgkey#1{\item{{\bf#1}.}}
%\newcount\qct\newcount\qcta\qct=0\qcta=1
%\def\pbgkey#1{\key{#1}\global\advance\qct
%       by\number\qcta
%       \immediate\write3
%       {\string\def\string#1\string{\number\qct\string}}}
\def\range{\operatorname{range}}
\def\rank{\operatorname{rank}}
\def\sign{\operatorname{sign}}
\def\trace{\operatorname{trace}}
\def\Pspan{\operatorname{span}}
%\immediate\openout 3=W01macro.tex
%%% Bibliographic references file

 \makeatletter
  \renewcommand{\theequation}{%
   \thesection.\arabic{equation}}
  \@addtoreset{equation}{section}
 \makeatother
%\topmatter
\title{Curvature tensors whose Jacobi or Szab\'o  operator is nilpotent
       on null vectors}
\begin{author}{Peter Gilkey\begin{thanks}{Research partially supported by the NSF (USA)
and the MPI (Leipzig)\newline{\it 2000 Mathematics Subject Classification. \rm Primary
53B20}\newline{\it Key words and phrases. \rm
   Algebraic curvature tensor,
          nilpotent Jacobi operator,\newline
   nilpotent Szab\'o operator, Lorentzian geometry, null vector}}\end{thanks}
\ and Iva Stavrov}
   \end{author}
%\address Mathematics Department, University of Oregon, Eugene Or
%97403 USA\endaddress
%\email gilkey\@darkwing.uoregon.edu, stavrov\@hopf.uoregon.edu\endemail
%\rightheadtext{Curvature Tensors}
%\leftheadtext{Gilkey and Stavrov}
\maketitle
\begin{abstract} We show that any $k$ Osserman Lorentzian
algebraic curvature tensor has constant sectional curvature and give an elementary proof
that any local $2$ point  homogeneous Lorentzian manifold has constant sectional
curvature. We also show that a Szab\'o Lorentzian covariant derivative algebraic curvature
tensor vanishes.\end{abstract}
\maketitle

%%%%%%%%%%%%%%%%%%%%%%%%%%%%%%%%%%%%%%%%%%%%%%%%%%%%%%%%%%%%%
\sethead\aref
\section{Introduction}\label{aref}
Let $V$ be a vector space
equipped with a symmetric inner product of signature $(p,q)$ and dimension $m=p+q\ge3$.
$V$ is said to be {\it Riemannian} if $p=0$ and {\it Lorentzian} if $p=1$. A
$4$ tensor $R\in\otimes^4V^*$ is said to be an {\it algebraic curvature tensor} if $R$
has the symmetries of the curvature of the Levi-Civita connection:
\begin{eqnarray*}
     &&R(x,y,z,w)=R(z,w,x,y)=-R(y,x,z,w),\text{ and }\\
     &&R(x,y,z,w)+R(y,z,x,w)+R(z,x,y,w)=0.\end{eqnarray*}
A $5$ tensor, which we denote symbolically by $\nabla R\in\otimes^5V^*$, is said to be a
{\it covariant derivative algebraic curvature tensor} if $\nabla R$ has the symmetries of
the covariant derivative of the curvature of the Levi-Civita connection:
\begin{eqnarray*}
&&\nabla  R(a,b,c,d;e)=-\nabla  R(b,a,c,d;e)=\nabla  R(c,d,a,b;e),\\
&&\nabla  R(a,b,c,d;e)+\nabla  R(a,c,d,b;e)+\nabla  R(a,d,b,c;e)=0,\text{ and}\\
&&\nabla  R(a,b,c,d;e)+\nabla  R(a,b,d,e;c)+\nabla  R(a,b,e,c;d)=0.\end{eqnarray*}
The Jacobi operator
${\Cal J}_R(x)$   and the Szab\'o operator ${\Cal S}_R(x)$ are the symmetric linear
operators on
$V$ defined by:
$$(\Cal{J}_R(x)y,w)=R(y,x,x,w)\text{ and }(\Cal{S}_{\nabla R}(x)y,w)=\nabla R(y,x,x,w;x).$$
 
In Section \ref{BREF}, we study the geometry of the Jacobi operator. 
Let $k$ be an index $1\le k\le m-1$ and let $\{e_1,...,e_k\}$ be an orthonormal basis for a
non-degenerate
$k$ dimensional subspace $\sigma\subset V$. The {\it higher order Jacobi operator} defined
by Stanilov and Videv \cite{refStanilovVidev} is the self-adjoint linear map of $V$ given
by:
$$\Cal{J}_R(\sigma):=\textstyle\sum_{1\le i\le k}(e_i,e_i)\Cal{J}_R(e_i);$$
this operator is independent of the particular orthonormal basis chosen for $\sigma$. The
algebraic curvature tensor
$R$ is said to be
$k$ Osserman if the eigenvalues of $\Cal{J}_R(\sigma)$ are constant on the Grassmannian of
non-degenerate $k$ planes in $V$. Note that if $R$ is $k$ Osserman, then $R$ is $m-k$
Osserman \cite{refGilkeyStanilovVidev}. Similarly, a pseudo-Riemannian manifold is said to
be $k$ Osserman if the eigenvalues of $\Cal{J}_R(\sigma)$ are constant on the Grassmannian
of non-degenerate $k$ planes in $TM$. It is conjectured
\cite{refOss} that a
$1$-Osserman Riemannian manifold is either flat or is locally a $2$ point homogeneous
space (i.e. a rank $1$ symmetric space). This is known if $m\equiv1$ mod $2$, if
$m\equiv2$ mod $4$, and if $m=4$ \cite{refChia}. Although there are some partial
additional results known, the general case remains open.

We say that linear map $A$ of $V$ is {\it nilpotent} if $A^m=0$
or equivalently if we have $\trace\{A^i\}=0$ for $1\le i\le m$. We complexify and extend
$g$,
$R$, and
$\nabla R$ to be complex multi-linear on $V_{\mathbb{C}}:=V\otimes\mathbb{C}$ so that we
can use analytic continuation. We say that a complex vector
$v$ is {\it null} if $(v,v)=0$; let $\Cal{N}$ be the set of all complex null vectors. In
Section \ref{BREF}, we prove the following two results:

\begin{theorem}\label{arefa} Let $R$ be an algebraic curvature tensor on a
vector space of arbitrary signature. If $R$ is $k$ Osserman, then $\Cal{J}_R(\cdot)$ is
nilpotent on $\Cal{N}$.\end{theorem}

\begin{theorem}\label{arefb} Let $R$ be an algebraic curvature tensor on a Lorentzian
vector space. If $\trace\{\Cal{J}_R(\cdot)^2\}=0$ on $\Cal{N}$, then $R$ has constant
sectional curvature.\end{theorem}

Theorems \ref{arefa} and \ref{arefb} imply the following result in the geometric setting
which was proved earlier \cite{refBokanBlazicGilkey,refGra} by different methods when
$k=1$. There is a similar result in the Riemannian setting if $2\le k\le m-2$
\cite{refGilkeya}. See also
\cite{refGilkeyStanilovVidev,refStanilovVidev}.

\begin{theorem}\label{arefc}
Let $(M,g)$ be a Lorentzian $k$ Osserman manifold. Then $(M,g)$
has constant sectional curvature.\end{theorem}

Let $S^\pm(V)$ and $S^\pm(M,g)$ be the pseudo-spheres and pseudo-sphere bundles of unit
timelike ($-$) and spacelike ($+$) vectors. We say that $(M,g)$ is a local $2$ point
homogeneous space if the local isometries of
$(M,g)$ act transitively on $S^\pm(M,g)$; this implies $\Cal{J}_R(\cdot)$ has
constant eigenvalues on $S^+(M,g)$ and on $S^-(M,g)$. The following result now follows
from Theorem \ref{arefc}:

\begin{theorem}\label{areff} If $(M,g)$ is a connected local 2 point homogeneous Lorentzian
manifold, then
$(M,g)$ has constant sectional curvature.
\end{theorem}

Szab\'o \cite{refSzabo} showed in
the Riemannian setting that if
$\Cal{S}_{\nabla R}(\cdot)$ has constant eigenvalues on $S(V)$, then $\nabla R=0$. He used this
observation to give an elementary proof that any local $2$ point homogeneous
Riemannian manifold is locally symmetric; of course more is true as $(M,g)$ is a local rank
$1$ symmetric space or is flat in this setting. This motivates the study of this operator
in the higher signature setting. In
Section \ref{CREF}, we study the geometry of the Szab\'o operator and prove the
following results:

\begin{theorem}\label{arefd}  Let $\nabla R$ be a covariant derivative algebraic curvature
tensor on a vector space of arbitrary signature. If the eigenvalues of $\Cal{S}_{\nabla R}$ are
constant on $S^+(V)$ and on $S^-(V)$, then $\Cal{S}_{\nabla R}(\cdot)$ is nilpotent
on $\Cal{N}$.\end{theorem}

\begin{theorem}\label{arefe} Let $\nabla R$ be a covariant derivative algebraic curvature
tensor on a Lorentzian vector space. If $\trace\{\Cal{S}_{\nabla R}(\cdot)^2\}$ is constant
on $S^\pm(V)$, then $\nabla R=0$.
\end{theorem}

We remark that one can use analytic continuation to show that if $p>0$ and $q>0$, then the
eigenvalues of $\Cal{S}_{\nabla R}$ are constant on $S^+(V)$ if and only the eigenvalues
of $\Cal{S}_{\nabla R}$ are constant on $S^-(V)$; we shall omit the proof in the interests
of brevity.

In both the Riemannian and the Lorentzian settings, if $\Cal{S}_{\nabla R}$ has constant
eigenvalues on $S^+(V)$ and on $S^-(V)$, then $\Cal{S}_{\nabla R}=0$. This can fail in the
higher signature setting. Note that if $\Cal{S}_{\nabla R}^2=0$, then $\Cal{S}_{\nabla R}$
has only the zero eigenvalue.

\begin{theorem}\label{arefg}  Let $V$ be a vector space of signature $(p,q)$, where
$p,q\ge2$. There exists a covariant derivative algebraic curvature tensor so
that $\Cal{S}_{\nabla R}^2(v)=0$ for all $v\in V$ and so that $\Cal{S}_{\nabla R}$ does not
vanish identically.
\end{theorem}

%In the Riemannian setting, topological methods are used to study when the Jacobi or
%Szab\'o operators have constant eigenvalues. We do not employ this machinery here. The
%situation in the higher signature setting is quite different and we refer to
%\cite{refGilkey} for a survey of what is known there.

\section{The Geometry of the Jacobi Operator\label{BREF}}

Let $\rho_R(x,y):=\trace\{z\rightarrow R(z,x)y\}$ be the Ricci tensor defined by $R$; we
then have $\rho_R(x,x)=\trace\{\Cal{J}_R(x)\}$.  We say $R$ is {\it Einstein} if there is a
constant $c_1$ so that $\rho_R(x,y)=c_1(x,y)$ for all $x,y\in V$. We adopt the notation of
\cite{refCarpenter} and say that $R$ is {\it $k$-stein} if there exist constants $c_i$ so
$\trace\{\Cal{J}_R(x)^i\}=c_i(x,x)^i$ for all $x\in V$. This definition is
motivated by the observation that $1$-stein and Einstein are equivalent notions. Note
that $R$ is $m$-stein if and only if $R$ is $1$-Osserman \cite{refGilkey}. We begin our
study of the geometry of the Jacobi operator with the following:

\begin{lemma}\label{brefa} Let $R$ be an algebraic curvature tensor on a vector space
of arbitrary signature.\begin{enumerate}
\smallskip\item If $R$ is $k$-stein, then
$\trace\{\Cal{J}_R(\cdot)^k\}=0$ on $\Cal{N}$.
\item If $R$ is $m$-stein, then
$\Cal{J}_R(\cdot)$ is nilpotent on $\Cal{N}$.
\item If $\trace\{\Cal{J}_R(\cdot)\}=0$ on $\Cal{N}$, then
$R$ is Einstein.
\end{enumerate}\end{lemma}

\medbreak\noindent{\bf Proof:} Let $R$ be $k$-stein. We use analytic continuation to see
the identity
$\trace\{\Cal{J}_R(x)^i\}=c_i(x,x)^i$ holds for complex vectors $x$ as well if $1\le i\le
k$. Assertions (1) and (2) then follow. Let
$x_1$ and $x_2$ be complex vectors so
$g(x_1,x_1)=1$, $g(x_2,x_2)=1$, and $g(x_1,x_2)=0$. Since the cross terms
cancel, we may expand
$${\Cal J}_R(x_1+\sqrt{-1}x_2)+\Cal{J}_R(x_1-\sqrt{-1}x_2)
  =2{\Cal J}_R(x_1)-2{\Cal J}_R(x_2).$$
Since $x_1\pm\sqrt{-1}x_2$ is a complex null vector, 
$\trace\{\Cal{J}_R(x_1\pm\sqrt{-1}x_2)\}=0$. This shows
$\trace\{\Cal{J}_R(x_1)\}=\trace\{\Cal{J}_R(x_2)\}$. It now follows
that $\rho_R(x,x)=c(x,x)$ for any complex vector
$x$ and consequently $R$ is Einstein.
\qedbox

\medbreak We remark that it is necessary to deal with complex null vectors in Lemma
\ref{brefa} to ensure that the statements are non-vacuous in the definite setting as there
are no real null vectors if $p=0$ or if $q=0$. It is not known if the
converse to assertion (2) holds, i.e. if $\Cal{J}_R(\cdot)$ is nilpotent on $\Cal{N}$,
then is $R$ is $m$-stein?

\medbreak\noindent{\bf Proof of Theorem \ref{arefa}} Let $R$ be $k$ Osserman and let
$x_1\in\Cal{N}$. We must show $\trace\{\Cal{J}_R(x_1)^k\}=0$ for all $k$. Since the inner
product on the complexification of $V$ is non-degenerate, we can choose
$x_2$ so $(x_1,x_2)\ne0$. Define linear functions $T_i$ by $T_i(x):=(x,x_i)$. Since
$T_1(x_1)=0$, $T_1(x_2)\ne0$, and $T_2(x_1)\ne0$,
the two linear functions $T_1$ and $T_2$ are linearly independent. We define:
$$W:=\ker(T_1)\cap\ker(T_2)=x_1^\perp\cap x_2^\perp.$$
If $a_1x_1+a_2x_2\in W$, then $0=(a_1x_1+a_2x_2,x_1)=a_2(x_2,x_1)$ so $a_2=0$ since
$(x_2,x_1)\ne0$. If $a_1x_1\in W$, then $0=(a_1x_1,x_2)=a_1(x_1,x_2)$ so $a_1=0$. Thus
$$\Pspan\{x_1,x_2\}\cap W=\{0\}.$$
Thus there is a basis $\{x_1,x_2,w_3,...,w_m\}$ for $V_{\mathbb{C}}$ so that the
subset
$\{w_3,...,w_m\}$ is a basis for $W$. Suppose that $w\in W$ and that $(w,w_i)=0$ for $3\le
i\le m$. Since $(w,x_1)=(w,x_2)=0$, this implies that $(w,v)=0$ for all $v$ and
hence $w=0$. This shows that the induced inner product on $W$ is non-degenerate.

As $k-1\le p+q-2=\dim W$, there is a non-degenerate $k-1$ plane
$\sigma\subset W$. Let 
$$x_t:=x_1+tx_2\text{ and }g(t):=(x_t,x_t)=(x_2,x_2)t^2+2t(x_1,x_2).$$
As $(x_1,x_2)\ne0$, $g(t)$ is a non-trivial polynomial of degree at most $2$.
Thus $g(t)$ has at most two roots; in particular, there exists $\varepsilon>0$ so
$g(t)\ne0$ for $t\in(0,\varepsilon)$. Let $\pi(t):=\sigma\oplus\Pspan\{x_t\}$;
$\pi(t)$ is a non-degenerate $k$ plane for
$t\in(0,\varepsilon)$. Let $t\in(0,\varepsilon)$.
As $R$ is $k$ Osserman, there
are universal constants $c_i$ so $c_i=\trace\{\Cal{J}_R(\tilde\pi)^i\}$ for
any non-degenerate real $k$ plane $\tilde\pi\subset V$. Again, analytic continuation
permits us to extend this relationship to the complex setting so:
$$c_i=\trace\{{\Cal J}_R(\pi(t))^i\}\text{ for }t\in(0,\varepsilon).$$
Since $\pi(t)=\sigma\oplus\Pspan\{x_t\}$ is an
orthogonal direct sum,
\begin{eqnarray*}
   \Cal{J}_R(\pi(t))&=&\Cal{J}_R(\sigma)+g(t)^{-1}\Cal{J}_R(x_t),\\
   c_i&=&\trace\{[{\Cal J}_R(\sigma)+g(t)^{-1}{\Cal J}_R(x_t)]^i\},\text{ and}\\
    g(t)^ic_i&=&g(t)^i\trace\{[{\Cal J}_R(\sigma)+g(t)^{-1}{\Cal J}_R(x_t)]^i\}\\
    &=&\trace\{[g(t){\Cal J}_R(\sigma)+{\Cal J}_R(x_t)]^i\}.\end{eqnarray*}
We take the limit as $t\downarrow0$
to complete the proof that $\trace\{{\Cal J}_R(x_1)^i\}=0$. \qedbox

\medbreak We have the following useful characterization
of $1$ Osserman algebraic curvature tensors in terms of the order of vanishing of
$\trace\{\Cal{J}_R(x)^k\}$ on the space of complex null vectors
$\Cal{N}$.

\begin{lemma}\label{brefb}  Let $R$ be an algebraic curvature tensor on a vector space of
arbitary signature. The following conditions are equivalent:
\begin{enumerate}
\item $R$ is $1$ Osserman
\item $\trace\{\Cal{J}_R(x+ty)^k\}=O(t^k)$ as $t\downarrow 0$ $\forall\
x\in\Cal{N}$,
$\forall\ y\in V_{\mathbb{C}}$, and
$\forall\ k$.
\end{enumerate}\end{lemma}

\medbreak{\bf Proof:} By definition, if $R$ is $1$ Osserman, then there exists a constant
$c_k$ so that
$\trace\{\Cal{J}_R(z)\}=c_k(z,z)^k$ for any $z\in V$. We complexify and use analytic
continuation to see this identity continues to hold for $z\in V_{\mathbb{C}}$.
We set $z=x+ty$ and $(x,x)=0$ to see
$$\trace\{\Cal{J}_R(x+ty)^k\}=c_k(x+ty,x+ty)^k=c_kt^k\{2(x,y)+t(y,y)\}^k.$$
thus assertion (1) implies assertion (2). Conversely, suppose that assertion (2) holds.
Let $e_1$ and $e_2$ be complex unit vectors.
We define a polynomial of degree $2k$
$f(t,s):=\trace\{\Cal{J}_R(te_1+se_2)^k\}$ in the
variables $(t,s)$. Since
$$f(t,s)=\trace\bigl\{\Cal{J}_R\{t(e_1\pm\sqrt{-1}e_2)+(s\mp t\sqrt{-1})e_2\}\bigr\}
    =O(s\mp\sqrt{-1}t)^k,$$
$(s\pm\sqrt{-1}t)^k$ divides $f(t,s)$. Thus
$f(t,s)=C(s^2+t^2)^k$. Since we may expand
\begin{eqnarray*}
f(t,s)&=&t^{2k}\trace\{\Cal{J}_R(e_1)^k\}+...+s^{2k}\trace\{\Cal{J}_R(e_2)^k\}\\
&=&C(s^2+t^2)^k=Ct^{2k}+...+Cs^{2k}\end{eqnarray*}
we have $\trace\{\Cal{J}_R(e_1)^k\}=\trace\{\Cal{J}_R(e_2)^k\}$ for all $k$ and all
complex unit vectors $e_i$. This shows that
$R$ is $k$ Osserman.
\qedbox

\medbreak\noindent{\bf Proof of Theorem \ref{arefb}} Let $V$ have signature $(1,q)$. We
suppose without loss of generality that
$q\ge2$. Let $\{e_0,...,e_q\}$ be an orthonormal basis for $V$, where $e_0$ is
timelike.  Let $\varepsilon_i:=(e_i,e_i)$; $\varepsilon_0=-1$ and $\varepsilon_i=+1$ for
$i\ge1$. We shall set
$R_{ijkl}:=R(e_i,e_j,e_k,e_l)$. By assumption, $\trace\{J_R(e_0\pm e_1)^2\}=0$. We have:
\begin{eqnarray}
   0&=&\textstyle\frac12\trace\{\Cal{J}_R(e_0+e_1)^2+\Cal{J}_R(e_0-e_1)^2\}\nonumber\\
   &=&\textstyle\sum_{ij}\varepsilon_i\varepsilon_j
     \{(R_{i00j}+R_{i11j})^2+(R_{i10j}+R_{i01j})^2\}\label{brefba}\end{eqnarray}
Let $j$ be an arbitrary index. The two terms in equation (\ref{brefba}) with $i=0$ and
$i=1$ are given by:
\begin{eqnarray*}
&&(R_{011j}+R_{000j})^2+(R_{001j}+R_{010j})^2=R_{011j}^2+R_{100j}^2\cr
   &=&(R_{111j}+R_{100j})^2+(R_{101j}+R_{110j})^2.\end{eqnarray*}
Since $\varepsilon_1\varepsilon_j+\varepsilon_0\varepsilon_j=0$,  the terms in equation
(\ref{brefba}) with $i=1$ and $i=0$ cancel. We may therefore restrict of the index $i$ in
equation to the range $2\le i\le q$. A similar argument shows that we may restrict $j$ to
the range $2\le j\le q$. Since $\varepsilon_i=\varepsilon_j=+1$ in this range, equation
(\ref{brefba}) shows that a sum of squares is zero. Consequently
$$R_{i11j}=-R_{i00j}\text{ and }R_{i10j}=-R_{j01i}\text{ for }2\le i,j\le q.$$
This holds for any orthonormal basis for $V$, where $e_0$ is timelike. It now follows
that $V$ has constant sectional curvature.
\qedbox

\section{The Geometry of the Szab\'o operator\label{CREF}}

\medbreak\noindent{\bf Proof of Theorem \ref{arefd}} Let $\nabla R$ be a covariant
derivative algebraic curvature tensor on a vector space of signature $(p,q)$. If $p>0$,
assume the eigenvalues of $\Cal{S}_{\nabla R}(\cdot)$ are constant on $S^-(V)$; the
argument is similar if $q>0$. If $k$ is odd, then $\trace\{\Cal{S}_{\nabla
R}(x)^k\}=(-1)^k\trace\{\Cal{S}_{\nabla R}(-x)^k\}$ and hence $\trace\{\Cal{S}_{\nabla
R}(x)^k\}=0$ on $S^-(V)$; analytic continuation then implies $\trace\{(\Cal{S}_{\nabla
R})^k\}$ vanishes identically. We may therefore suppose $k$ even. We rescale to see there
are constants $c_j$ so
\begin{equation}\label{crefaa}
\trace\{\Cal{S}_{\nabla R}(x)^{2j}\}-c_j(x,x)^{3j}=0\end{equation}
on the open subset of all timelike vectors. Analytic continuation then implies equation
(\ref{crefaa}) holds for all vectors; we take $x\in\Cal{N}$ to complete the proof.
\qedbox

\begin{remark}\label{crefh} The same argument used to prove Lemma \ref{brefb} extends to
show that $\nabla R$ is Szab\'o if and only if we have that $\trace\{S_{\nabla
R}(x+ty)^{2k}\}=O(t^{3k})$ and that $\trace\{S_{\nabla R}(x+ty)^{2k-1}\}=0$ for
every $x\in\Cal{N}$.
\end{remark}

We  begin the proof of Theorem \ref{arefe} with a technical result.
\begin{lemma}\label{crefe} Let $\nabla R$ be a covariant derivative algebraic curvature
tensor on a Lorentzian vector space. If $\trace\{\Cal{S}_{\nabla R}(\cdot)^2\}$ is
constant on $S^-(V)$, $\Cal{S}_{\nabla R}=0$.\end{lemma}

\medbreak\noindent{\bf Proof:} As in the proof of Theorem \ref{arefb}, let
$\Cal{B}:=\{e_0,e_1,...,e_q\}$ be an orthonormal basis for $V$, where $e_0$ is timelike
and $e_i$ is spacelike for $i>0$. Let
$\theta$ be a real parameter. We define a new orthonormal basis $\Cal{B}(\theta)$ by:
\begin{eqnarray*}
&&e_0(\theta):=\cosh\theta\cdot e_0+\sinh\theta\cdot e_1,\ 
 e_1(\theta):=\sinh\theta\cdot e_0+\cosh\theta\cdot e_1,\\
&&e_i(\theta):=e_i\text{ for }i\ge2.\end{eqnarray*}
We have a constant $C$ so that:
\begin{eqnarray}
C&=&\trace\{\Cal{S}_{\nabla R}(e_0(\theta))^2\}\nonumber\\&=&
  \textstyle\sum_{1\le i,j\le q}\nabla R(e_i(\theta),e_0(\theta),e_0(\theta),e_j(\theta);
e_0(\theta))^2.\label{crefea}\end{eqnarray}
As $\cosh\theta=\frac12(e^\theta+e^{-\theta})$
and $\sinh\theta=\frac12(e^{\theta}-e^{-\theta})$, we may expand:
\begin{eqnarray*}
&&\nabla R(e_i(\theta),e_0(\theta),e_0(\theta),e_j(\theta);e_0(\theta))
=\textstyle\sum_{-5\le\nu\le5}a_{ij,\nu}e^{\nu\theta},\\
   &&C=\textstyle\sum_{1\le i,j\le q}\{a_{ij,5}\}^2e^{10\theta}+O(e^{9\theta}),\\
   &&Ce^{-10\theta}=\textstyle\sum_{1\le i,j\le q}\{a_{ij,5}\}^2+O(e^{-\theta}).
\end{eqnarray*}
We take the limit as $\theta\rightarrow\infty$ to see
$\textstyle\sum_{ij}\{a_{ij,5}\}^2=0$ and hence $a_{ij,5}=0$ for all $i,j$; similarly
$a_{ij,-5}=0$. Similarly we have
$a_{ij,\nu}=0$ for $\nu\ne0$. Consequently
\begin{equation}\label{crefec}
\nabla R(e_2,e_0(\theta),e_0(\theta),e_2;e_0(\theta))=a_{22,0}
\end{equation}
is independent of $\theta$. On the other hand, since there are three terms involving $\theta$, the powers of
$e^\theta$ which appear in this expression are odd. Thus $a_{22,0}=0$ so
$$\nabla R(e_2,e_0,e_0,e_2;e_0)=0.$$
Similarly we conclude
$\nabla R(e_i,e_0,e_0,e_i;e_0)=0$ for any $i\ge1$.
We polarize to see $\nabla R(e_i,e_0,e_0,e_j;e_0)=0$ for any $i,j$; the vanishing being
automatic if $i=0$ or $j=0$. Thus
$\Cal{S}_{\nabla R}(e_0)=0$. As $e_0$ was arbitary, $\Cal{S}_{\nabla R}(\cdot)=0$ on
$S^-(V)$. Rescaling and analytic continuation then imply
$\Cal{S}_{\nabla R}(\cdot)=0$ on $V$.
\qedbox

\medbreak We complete the proof of Theorem \ref{arefe} by proving:
\begin{lemma}\label{creff} Let $\nabla R$ be a covariant derivative algebraic curvature
tensor on a vector space of arbitrary signature. If $\Cal{S}_{\nabla R}=0$, then $\nabla
R=0$.\end{lemma}

\medbreak\noindent{\bf Proof:} Since we are in the purely algebraic setting, we only have
the pointwise vanishing that
$\nabla R=0$. Thus we can not appeal to an argument using Jacobi fields such as that
given in Besse \cite{refBesse}. Instead, we use
an argument based on the curvature symmetries given above; see, for example Vanhecke and
Willmore \cite{refVW} in the Riemannian setting. We
polarize the identity
\begin{equation}\label{creffa}
\nabla R(x,y,y,w;y)=0\text{ for all }w,x,y.\end{equation}
by setting $y(t):=y+tx$ and expanding in terms of powers of $t$. We set the term which is
linear in $t$ to zero and then use the curvature symmetries to compute:
\begin{eqnarray}
0&=&\nabla R(x,x,y,w;y)+\nabla R(x,y,x,w;y)+\nabla R(x,y,y,w;x)\nonumber\\
 &=&0+\nabla R(x,y,x,w;y)-\nabla R(x,y,w,x;y)-\nabla R(x,y,x,y;w)\nonumber\\
&=&-2\nabla R(x,y,w,x;y)+\nabla R(x,y,y,x;w).\label{creffb}\end{eqnarray}
Set $w=x$ in equation (\ref{creffa})  to see:
\begin{equation}\label{creffc}
\nabla R(x,y,y,x;y)=0\text{ for all }x,y.\end{equation}
Set $y(t):=y+tw$ and expand equation
(\ref{creffc}) in terms of powers of $t$; the term which is linear in $t$ then
yields the identity:
\begin{equation}\label{creffD}
0=2\nabla R(x,y,w,x;y)+\nabla R(x,y,y,x;w)\text{ for all }w,x,y.\end{equation}
We add equations (\ref{creffb}) and (\ref{creffD}) to see
\begin{equation}\label{creffe}
0=2\nabla R(x,y,y,x;w)\text{ for all }x,y,w.\end{equation}
Set $y(t):=y+tz$ and expand equation (\ref{creffe}) in terms of 
powers of $t$. The term which is linear in $t$ then yields the identity:
\begin{equation}\label{creffg}
0=\nabla R(y,x,x,z;w)\text{ for all }x,y,z,w\in V.\end{equation}
Set $x(t)=x+tv$ and expand equation (\ref{creffg}) in terms of powers of $t$. The
term which is linear in $t$ then yields the identity:
\begin{equation}\label{creffh}
0=\nabla R(y,x,v,z;w)+\nabla R(y,v,x,z;w)\text{ for all }v,w,x,y,z\in V.
\end{equation}
We use equation (\ref{creffh}) and the curvature symmetries to complete the proof:
\begin{eqnarray*}
   0&=&\nabla R(y,x,v,z;w)+\nabla R(y,v,z,x;w)+\nabla R(y,z,x,v;w)\\
    &=&\nabla R(y,x,v,z;w)-\nabla R(y,v,x,z;w)+\nabla R(y,z,x,v;w)\\
    &=&\nabla R(y,x,v,z;w)+\nabla R(y,x,v,z;w)-\nabla R(y,x,z,v;w)\\
    &=&3\nabla R(y,x,v,z;w)\text{ for all }v,w,x,y,z\in V.\ \qedbox\end{eqnarray*}

\medbreak\noindent{\bf Proof of Theorem \ref{arefg}} Let $\tilde L$ be a completely
symmetric trilinear form and
$L$ be a symmetric bilinear form on $V$. We define a $5$ tensor:
\begin{eqnarray}
\tilde R(x,y,z,w;v):&=&\tilde L(v,y,z)L(x,w)-\tilde L(v,x,z)L(y,w)\cr
                       &+&\tilde L(v,x,w)L(y,z)-\tilde L(v,y,w)L(x,z).
\label{crefga}\end{eqnarray}
It is immediate that $\tilde R(x,y,z,w;v)=-\tilde R(y,x,z,w;v)$. We interchange the roles
of $x$ and $z$ and of $y$ and $w$ to check:
\begin{eqnarray*}
\tilde R(z,w,x,y;v)&=&\tilde L(v,w,x)L(z,y)-\tilde L(v,z,x)L(w,y)\cr
                       &+&\tilde L(v,z,y)L(w,x)-\tilde L(v,w,y)L(z,x)\cr
&=&\tilde R(x,y,z,w;v).\end{eqnarray*}
We check the first Bianchi identity is
satisfied by computing:
\begin{eqnarray*}
&&\tilde R(x,y,z,w;v)+R(x,z,w,y;v)+R(x,w,y,z;v)\\
&=&\tilde L(v,y,z)L(x,w)+\tilde L(v,z,w)L(x,y)+\tilde L(v,w,y)L(x,z)\\
 &-&\tilde L(v,x,z)L(y,w)-\tilde L(v,x,w)L(z,y)-\tilde L(v,x,y)L(z,w)\\
 &+&\tilde L(v,x,w)L(y,z)+\tilde L(v,x,y)L(w,z)+\tilde L(v,x,z)L(y,w)\\
 &-&\tilde L(v,y,w)L(x,z)-\tilde L(v,z,y)L(x,w)-\tilde L(v,w,z)L(x,y)\\&=&0.
\end{eqnarray*}
We check the second Bianchi identity is satisfied by computing:
\begin{eqnarray*}
&&\tilde R(x,y,z,w;v)+\tilde R(x,y,w,v;z)+\tilde R(x,y,v,z;w)\\
&=&\tilde L(v,y,z)L(x,w)+\tilde L(z,y,w)L(x,v)+\tilde L(w,y,v)L(x,z)\\
&-&\tilde L(v,x,z)L(y,w)-\tilde L(z,x,w)L(y,v)-\tilde L(w,x,v)L(y,z)\\
&+&\tilde L(v,x,w)L(y,z)+\tilde L(z,x,v)L(y,w)+\tilde L(w,x,z)L(y,v)\\
&-&\tilde L(v,y,w)L(x,z)-\tilde L(z,y,v)L(x,w)-\tilde L(w,y,z)L(x,v)\\&=&0.
\end{eqnarray*}
Thus $\tilde R$ is an covariant derivative algebraic curvature tensor.
Furthermore,
\begin{eqnarray*}
  (\Cal{S}_{\nabla R}(x)y,w):&=&\tilde R(y,x,x,w;x)\\
   &=&\tilde L(x,x,x)L(y,w)-\tilde L(x,y,x)L(x,w)\\
                       &+&\tilde L(x,y,w)L(x,x)-\tilde L(x,x,w)L(y,x).\end{eqnarray*}
Let $\tilde\phi_x$ and $\phi$ be the associated self-adjoint operators;
$\tilde L(x,y,z)=(\phi_x(y),z)$ and
$L(x,y)=(\phi(x),y)$. We then have

\begin{eqnarray*}
\Cal{S}_{\nabla R}(x)y&=&\tilde L(x,x,x)\phi(y)-\tilde L(y,x,x)\phi(x)\\
   &+&L(x,x)\tilde \phi_x(y)-L(y,x)\tilde\phi_x(x).\end{eqnarray*}
We suppose $V$ has signature $(p,q)$. Let $\{e_i^\pm\}$ be an orthonormal basis for $V$
so $\Pspan\{e_1^-,...,e_p^-\}$ is a maximal timelike subspace and so
$\Pspan\{e_1^+,...,e_q^+\}$ is an orthogonal maximal spacelike subspace of $V$. Let
$\varepsilon$, $\delta$, and $\varrho$ be choices of $\pm1$ signs. For $1\le i,j,k\le 2$,
we define:
\begin{eqnarray*}
&&L(e_i^\varepsilon,e_j^\delta)=\delta_{ij},\qquad\qquad\quad\phantom{.}
\tilde L(e_i^\varepsilon,e_j^\delta,e_k^\varrho)=\delta_{ijk},\\
&&\phi(e_i^\pm)=-e_i^-+e_i^+,\text{ and }\quad
  \tilde\phi_{e_i^\pm}(e_j^\pm)=\delta_{ij}\phi(e_j^\pm).\end{eqnarray*}
We set $\tilde L(\cdot)=0$ and $L(\cdot)=0$ if any index is greater than $2$. It is
immediate from the definition that
$\Cal{S}_{\nabla R}(x)^2=0$. On the other hand
$$(\Cal{S}_{\nabla R}(e_1^\pm)e_2^\pm,e_2^+)\ne0\text{ so }\Cal{S}_{\nabla R}(x)\ne0.\
\qedbox$$

%\sethead\DREFf
%\head{\DREFf\ Conclusion}\endhead 
%We showed that the condition
%$\Cal{J}_R(\cdot)$ is nilpotent on $\Cal{N}$ arises naturally when
%considering $k$ Osserman algebraic curvature tensors and
%implies constant sectional curvature in the Lorentzian setting; we also showed that local
%$2$ point homogeneous spaces have constant sectional curvature in the Lorentzian setting. 
%Similarly we showed that the condition $\Cal{S}_{\nabla R}(\cdot)$ is nilpotent on
%$\Cal{N}$ arises naturally in the study of covariant derivative algebraic curvature
%tensors whose Szab\'o operator has constant eigenvalues.

\medbreak
\par Department of Mathematics
\par University of Oregon
\par Eugene, Or 97403 USA
\par gilkey@darkwing.uoregon.edu
\par stavrov@hopf.uoregon.edu

\enddocument